\documentclass[]{article}

\usepackage[english]{babel}
\usepackage[T1]{fontenc}
\usepackage[utf8]{inputenc}
\usepackage[a4paper, margin=3cm]{geometry}
\usepackage{lmodern}
\usepackage{amsmath}
\usepackage{amsfonts}
\usepackage{dsfont}
\usepackage{amsthm}
\usepackage{mathtools}
\usepackage[authoryear]{natbib}
\usepackage[ruled,linesnumbered]{algorithm2e}
\usepackage{todonotes}

\numberwithin{equation}{section}

\newcommand{\DEF}{\coloneqq}
\newcommand{\R}{\mathbb{R}}
\newcommand{\N}{\mathbb{N}}
\newcommand{\K}{\mathcal{K}}
\newcommand{\SK}{\mathcal{S}_{q}(\mathcal{K})}

\begin{document}
	
	
\title{A Multigrid Preconditioner for Tensor Product Spline Smoothing}
\author{Martin Siebenborn\thanks{University of Hamburg, Department of Mathematics, Bundesstraße~55, 20146 Hamburg, Germany, Email: martin.siebenborn@uni-hamburg.de} \and Julian Wagner\thanks{Trier University, DFG-RTG Algorithmic Optimization, Universitätsring 15, 52496 Trier, Germany, Email: wagnerj@uni-trier.de} 
}
\date{}
\maketitle


\begin{abstract}
\noindent
Uni- and bivariate data smoothing with spline functions is a well established method in nonparametric regression analysis.
The extension to multivariate data is straightforward, but suffers from exponentially increasing memory and computational complexity.
Therefore, we consider a matrix-free implementation of a geometric multigrid preconditioned conjugate gradient method for the regularized least squares problem resulting from tensor product B-spline smoothing with multivariate and scattered data.
The algorithm requires a moderate amount of memory and is therefore applicable also for high-dimensional data.
Moreover, for arbitrary but fixed dimension, we achieve grid independent convergence which is fundamental to achieve algorithmic scalability.
\end{abstract}




\section{Introduction}
\label{sec:introduction}


In many statistical applications it is fundamental to investigate the relationship between explanatory variables and a variable of interest, which is generally modeled by a function of the covariates.
This function attempts to capture the important patterns in the data while leaving out noise and other insignificant structures.
This concept is known under several names like smoothing, (nonparametric) regression, or surface fitting.
The input data is often multivariate and scattered, i.e.\ there is no inherent structure.
To represent complex, e.g.\ multivariate and highly nonlinear data, the modeled function has to be very flexible in order to allow data-driven estimation of the complex effects.
In one and two dimensions, there exist efficient smoothing methods, based on spline functions.
We refer the reader to, e.g., \cite{Eilers1996}, \cite{Ruppert2003}, \cite{Wand2008}, \cite{Fahrmeir2013}.
Unfortunately, the straightforward extension of these spline-based methods to multiple input variables suffers from an exponential growth of the number of parameters to be estimated within the spatial dimension.
This issue is often referred to as the curse of dimensionality (cf. \cite{Bellman1957}).
Therefore, the computational and especially memory complexity of the related estimation procedure becomes unjustifiably large, even for moderate spatial dimensions, such that smoothing methods are rarely applied for covariates of dimension larger than two (cf. \cite{Fahrmeir2013}, p. 531).
For gridded covariates an efficient algorithm, that extends the spline smoothing method of \cite{Eilers1996}, is implemented in \cite{Eilers2006}.
However, for scattered covariates there is no satisfactory approach.
The challenge is to deal with the large-scale linear systems which arise from the spline smoothing for scattered data sets with increasing covariate dimensions.

In order to overcome this issue, we apply two techniques in the solution process and investigate their performance.
First, the matrix corresponding to the linear system is never explicitly formed and stored in memory.
When using an iterative solution algorithm like the conjugate gradient (CG) method only the matrix-vector product is required, but not matrix entries explicitly.
A second important ingredient is a suitable preconditioner for the linear system.
For this purpose, we introduce multigrid techniques to the spline smoothing problem.
With increasing space dimension and grid resolution the performance of the CG method usually deteriorates.
We thus investigate a geometric multigrid preconditioned and matrix-free CG iteration, which significantly reduces memory and computational complexity compared to common estimation methods.
Furthermore, the grid-independent convergence of the geometric multigrid is mandatory in order to achieve algorithmic scalability in the sense that simultaneously doubling the degrees of freedom and the number of processors leads to constant algorithmic running times.

The remainder of this paper is organized as follows.
In Section \ref{sec:regression}, we formulate the spline smoothing problem in multiple dimensions based on tensor product B-splines and formulated the underlying large-scale linear system.
In Section \ref{sec:multigrid}, we develop a matrix-free implementation of a multigrid preconditioned conjugate gradient method to solve the large-scale system with comparatively small memory and computational complexity.
In Section \ref{sec:numerics}, we apply the proposed method, also compared to traditional methods, on a test data set.
Finally, Section \ref{sec:conclusion} gives a conclusion of the paper.


\section{Spline Smoothing in Multiple Dimensions}
\label{sec:regression}


In statistics, smoothing a given data set describes the process of constructing a function that captures the important patterns in the data while leaving out noise and other fine scaled structures.
More precisely, for given data $\lbrace (x_i,y_i) \in \R^P \times \R \,|\, i=1,\ldots,n \rbrace$, where the $y_i \in \R$ are observations of a continuous response variable and the $x_i \in \R^P$ represent the corresponding value of a continuous covariate, we seek a smooth, but further unspecified function $s \colon \Omega \subset \R^P \rightarrow \R$ such that
\begin{align*}
	y_i = s(x_i) + \varepsilon_i, \ i=1,\ldots,n .
\end{align*}
A common assumption is that $\varepsilon_1,\ldots,\varepsilon_n$ are independent and identically distributed (i.i.d.) random errors with zero mean, common variance $\sigma_{\varepsilon}^2$ and assumed to be independent of the covariates.
We begin by defining a spline basis which forms the underlying, linear space for the representation of $s$.


\subsection{Tensor Product Splines}


Let $\Omega \DEF [a,b]$ be a bounded and closed interval partitioned by the knots
\begin{align*}
	\K \DEF \lbrace a=\kappa_0 < \ldots < \kappa_{m+1}=b \rbrace .
\end{align*}
Let $\mathcal{C}^{q}(\Omega)$ denote the space of $q$-times continuously differentiable functions and let $\mathcal{P}_q(\Omega)$ denote the space of polynomials of degree $q$.
We call the function space
\begin{align*}
	\SK \DEF \lbrace s \in \mathcal{C}^{q-1}(\Omega) :  s|_{\left[ \kappa_{j-1},\kappa_{j} \right]} \in \mathcal{P}_q \left( \left[ \kappa_{j-1},\kappa_{j} \right] \right), \ j=1,\ldots,m+1 \rbrace
\end{align*}
the space of spline functions of degree $q \in \N_0$ with knots $\K$.
It is a finite dimensional, linear space of dimension $J \DEF dim \left( \SK \right) = m+q+1$.
With 
\begin{align*}
	\lbrace \varphi_{j,q} : j=1,\ldots,J \rbrace 
\end{align*}
we denote its B-spline basis (cf. \cite{Boor1978}), which is well suited for numerical applications.
To extend the spline concept to $P$-dimensional covariates, a tensor product approach is common.
Let $\mathcal{S}_{q_p}(\K_p)$ be the spline space for the $p$-th covariate, $p=1,\ldots,P$, and let 
\begin{align*}
	\lbrace \varphi_{j_p,q_p}^p : j_p=1,\ldots,J_p \rbrace 
\end{align*}
denote its B-spline basis.
The function 
\begin{align*}
	\varphi_{j,q} \colon \Omega \DEF \Omega_1 \times \ldots \times \Omega_P \rightarrow \R , \ \varphi_{j,q}(x) = \prod_{p=1}^{P} \varphi^p_{j_p,q_p}(x^p) ,
\end{align*}
where $j \DEF (j_1,\ldots,j_P){'}$ and $q \DEF (q_1,\ldots,q_P){'}$ are multiindices, is called tensor product B-spline.
We define the space of tensor product splines as their linear combination
\begin{align} \label{Equation:tensor_product_splines}
	\SK \DEF \text{span} \lbrace \varphi_{j,q} : 1 \leq j \leq J \DEF (J_1, \ldots, J_P){'} \rbrace ,
\end{align}
which is then a $K \DEF \prod_{p=1}^{P} J_p$ dimensional linear space.
Note, that we use the same symbol for the univariable and the tensor product spline space as well as for B-splines and tensor product B-splines.
The difference is that we make use of the multiindex notation for the tensor products.
Every tensor product spline $s \in \SK$ therefore has a unique representation
\begin{align*}
	s = \sum\limits_{1 \leq j \leq J} \alpha_j \varphi_{j,q}
\end{align*} 
and for computational reasons we uniquely identify the set of multiindices $\lbrace j \in \N^P : 1 \leq j \leq J \rbrace$ in descending lexicographical order as $\lbrace 1,\ldots,K \DEF \prod_{p=1}^{P} J_p \rbrace$ such that
\begin{align} \label{Equation:BSplineRepresentation}
	s = \sum\limits_{1 \leq j \leq J} \alpha_j \varphi_{j,q} = \sum\limits_{k=1}^{K} \alpha_k \varphi_{k,q} .
\end{align}
We are now prepared to formulate the so called smoothing spline problem.


\subsection{Tensor Product Smoothing Spline}


Spline smoothing, also known as regularized or penalized spline regression, is a popular method in data smoothing.
We define (tensor product) smoothing splines in $\SK$ as solution of
\begin{align} \label{Equation:SmoothingSpline}
	\min \limits_{s \in \SK} \sum\limits_{i=1}^{n} \left( s(x_i)-y_i \right)^2 + \lambda \int\limits_{\R^P} \left( \sum\limits_{p_1=1}^{P} \sum\limits_{p_2=1}^{P} \dfrac{\partial^2}{\partial x_{p_1} \partial x_{p_2}} s(x) \right)^2 \mathrm{d}x .
\end{align}
The objective function consists of several parts.
One the one hand, we have a least squares fitting term
\begin{align*} 
	\mathcal{LS}(s) \DEF \sum\limits_{i=1}^{n} \left( s(x_i)-y_i \right)^2
\end{align*}
that measures the goodness-of-fit of the smoothing spline to the given observations.
On the other hand, we have a regularization term
\begin{align*}
	\mathcal{R}(s) \DEF \int\limits_{\R^P} \left( \sum\limits_{p_1=1}^{P} \sum\limits_{p_2=1}^{P} \dfrac{\partial^2}{\partial x_{p_1} \partial x_{p_2}} s(x) \right)^2 \mathrm{d}x ,
\end{align*}
that penalizes the roughness of the spline function.
The term $\mathcal{R}(s)$ is weighted by a regularization or smoothing parameter $\lambda>0$, which balances the two competitive terms $\mathcal{LS}(s)$ and $\mathcal{R}(s)$.
For $\lambda \rightarrow 0$ the smoothing spline tends to be the common least squares spline, which heavily overfits, or even interpolates, the observations.
Conversely, for $\lambda \rightarrow \infty$, to much impact is given to the regularization term such that the smoothing spline tends to be the least squares hyperplane (cf. \cite{Green1993}, p. 159).
More details can be found in the monographs \cite{Eubank1988}, \cite{Wahba1990}, and \cite{Green1993}.

Using the unique B-spline representation \eqref{Equation:BSplineRepresentation} we obtain
\begin{align*}
	\mathcal{LS}(s) = \Vert \Phi\alpha - y \Vert_2^2 ,
\end{align*}
where $\Phi \in \R^{n \times K}$ is element wise defined as $\Phi[i,k] \DEF \varphi_{k,q}(x_i)$ and 
\begin{align*}
	\mathcal{R}(s) = \alpha{'} \left( \sum\limits_{r \in \N_0^P ; \vert r \vert = 2} \frac{2}{r!} \Psi_r \right) \alpha ,
\end{align*}
where each $\Psi_{r} \in \R^{K \times K}$ is element wise defined as $\Psi_{r}[k,\ell] = \left\langle \partial^{r} \varphi_{k,q} , \partial^{r} \varphi_{\ell,q} \right\rangle_{L^2(\Omega)}$.
Defining
\begin{align*}
	\Lambda \DEF \sum\limits_{r \in \N_0^P ; \vert r \vert = 2} \frac{2}{r!} \Psi_r \in \R^{K \times K}
\end{align*}
the smoothing spline \eqref{Equation:SmoothingSpline} is equivalently formulated in terms of the B-spline coefficients as
\begin{align*}
	\min\limits_{\alpha \in \R^K} \Vert \Phi\alpha-y\Vert_2^2 + \lambda \alpha{'}\Lambda\alpha .
\end{align*} 
This is a regularized least squares problem and its solution is given by the solution of the linear system
\begin{align} \label{Equation:SmoothingSplineLinearSystem}
	\left( \Phi{'} \Phi + \lambda \Lambda \right)\alpha  \stackrel{!}{=} \Phi{'}y .
\end{align}
The main focus of this paper is a memory efficient evaluation and solution of \eqref{Equation:SmoothingSplineLinearSystem}.


\subsection{Curse of Dimensionality}


Spline smoothing is known to suffer from the so called curse of dimensionality, which describes an exponential growth of the number of B-spline coefficients $K$ within the dimension of the covariates $P$.
Choosing $J_p=35$ B-spline functions for each direction $p=1,\ldots,P$, which is a quite realistic number, the resulting tensor product B-spline is given by
\begin{align*}
	K = \prod_{p=1}^{P} J_p = 35^P
\end{align*}
parameters.
Table \ref{Table:CurseOfDimensionality} gives an overview on the increasing number of degrees of freedom under increasing dimensionality with which we are dealing with in this work.
\begin{table}[htb]
	\centering
	\begin{tabular}{c|c|c|c|c}
		& $P=2$ & $P=3$ & $P=4$ & $P=5$\\
		\hline
		$K$ & $1{.}225$ & $42{.}875$ & $1{.}500{.}625$ & $52{.}521{.}875$ \\
	\end{tabular}
	\caption{Number of B-spline coefficients for varying spatial dimension $P$.}
	\label{Table:CurseOfDimensionality}
\end{table}
Obviously, the method of smoothing splines becomes impracticable for dimensions $P \geq 3$ (cf. \cite{Fahrmeir2013}, p. 531), since solving the large-scale linear system \eqref{Equation:SmoothingSplineLinearSystem} requires an unjustified computational effort.
One reason for this is that, even when the coefficient matrix $\Phi{'}\Phi+\lambda\Lambda$ is stored in a sparse format, the exponential growth ensures that the memory capacity of a common digital computer is exceeded already for moderate spatial dimensions $P$.

Therefore, to make tensor product spline smoothing practicable for increasing spatial dimensions $P$, we need to develop computational and memory efficient methods to solve the large-scale linear system \eqref{Equation:SmoothingSplineLinearSystem}.


\subsection{Tensor Product Properties}


Before focusing on efficient solution methods in the next section, we state some important properties on the occurring matrices that are fundamental for the proposed algorithm and are based on the tensor product nature of the underlying splines.
First, we define for each spatial direction $p=1,\ldots,P$ the matrix
\begin{align*}
	\Phi_p \in \R^{n \times J_p} , \ \Phi_p[i,j_p] \DEF \varphi^p_{j_p,q_p}(x_i^p) ,
\end{align*}
which corresponds to the matrix $\Phi$ with a one-dimensional covariate in direction $p$.
Then, because of the scattered data structure, it holds
\begin{align*}
	\Phi{'} = \bigodot\limits_{p=1}^P \Phi_p{'} ,
\end{align*}
where $\odot$ denotes the Khatri-Rao product.
For each Gramian matrix $\Psi_r$ we analogously define a Gramian matrix for each spatial direction of the respective derivatives, that is
\begin{align*}
	\Psi_{r_p}^{p} \in \R^{J_p \times J_p} \ ,\Psi_{r_p}^{p}[j_p,\ell_p]= \left\langle \partial^{r_p} \varphi_{j_p,q_p}^p , \partial^{r_p} \varphi_{\ell_p,q_p}^p \right\rangle_{L^2(\Omega_p)} .
\end{align*}
Then it holds
\begin{align*}
	\Psi_{r} = \bigotimes\limits_{p=1}^{P} \Psi^{p}_{r_p}
\end{align*}
due to the tensor property.

A further important property of uniform B-splines in one variable is the subdivision formula.
Let $\mathcal{S}_q(\K^{2h})$ and $\mathcal{S}_q(\K^{h})$ denote univariable spline spaces with uniform knot set $\K^{2h}$ and $\K^{h}$ with mesh sizes $2h$ and $h$, respectively.
Then it holds (cf. \cite{Hoellig2003}, p. 32) 
\begin{align} \label{Equation:SubdivisionFormula}
	\varphi_{j,q}^{2h} = \frac{1}{2^q} \sum\limits_{i=0}^{q+1} \binom{q+1}{i} \varphi_{2j-(q+1)+i,q}^{h}
\end{align}
and consequently for a uniform spline $s \in \mathcal{S}_q(\K^{2h})$ we obtain
\begin{align*}
	s = \sum\limits_{j=1}^{J^{2h}} \alpha_{j}^{2h} \varphi_{j,q}^{2h} = \sum\limits_{j=1}^{J^{2h}} \alpha_{j}^{2h} \frac{1}{2^q} \sum\limits_{i=0}^{q+1} \binom{q+1}{i} \varphi_{2j-(q+1)+i,q}^{h} .
\end{align*}
Since $s \in \mathcal{S}_q(\K^{h})$, it holds
\begin{align*}
	s = \sum\limits_{j=1}^{J^{h}} \alpha_{j}^{h} \varphi_{j,q}^{h}
\end{align*}
and the B-spline coefficients for the different meshes are therefore related through $\alpha^{(h)} = I_{2h}^{h} \alpha^{2h}$, where $I_{2h}^{h} \in \R^{J^{h} \times J^{2h}}$ is element wise defined as
\begin{align*}
	I_{2h}^{h} [i,j] \DEF \frac{1}{2^q} \binom{q+1}{i-2j+q+1}, \ i=1,\ldots,J^{h}, \ j= 1,\ldots,J^{2h}.
\end{align*}
Because of the tensor property of B-splines, the formula carries over to multivariable B-splines and the corresponding B-spline coefficients are related through the matrix
\begin{align} \label{Equation:SubdivisionMatrixMulti}
	I_{2h}^{h} = \bigotimes_{p=1}^{P} I_{2h_p}^{h_p} \in \R^{ K^{h} \times K^{2h} },
\end{align}
where $h \DEF (h_1,\ldots,h_P){'}$ denotes the mesh vector.


\section{Geometric Multigrid Preconditioner}
\label{sec:multigrid}


In this section we focus on the numerical behavior of the normal equation (\ref{Equation:SmoothingSplineLinearSystem}), which is equivalent to solving the smoothing spline problem (\ref{Equation:SmoothingSpline}).
Since the coefficient matrix $A \DEF \Phi{'}\Phi + \lambda \Lambda$ is symmetric and positive definite, the conjugated gradient (CG) method is the straightforward choice.
For our application it is particularly important, that for the CG algorithm only matrix-vector products are required but not explicit entries of the matrix.
From a computational point of view assembling $A$ in a sparse format is infeasible for the targeted problem sizes.
This limits the choice for linear solvers and preconditioners.
It is common practice to apply a preconditioner to the CG method in order to speed up convergence.
Since the convergence rate depends on the condition of the matrix we encounter a worsening while improving the approximation quality by refining the spline space \eqref{Equation:tensor_product_splines}.
This behavior is more dramatic, the higher the spatial dimension is, which is reflected in the numerical test cases in Section \ref{sec:numerics}.
For the choice of a preconditioner it is again important, that only matrix-vector products are involved.
This cancels out the well-established, incomplete Cholesky factorizations of $A$ as a preconditioner and we thus concentrate on geometric multigrid methods in this paper.

The origins of multigrid methods date back into the late 70th (see for instance \cite{Brandt1977}, \cite{Hackbusch1978}, or \cite{Trottenberg2000} for an overview).
Since then, they are successfully applied in the field of partial differential equations.
The main idea is to build a hierarchy of - in a geometrical sense - increasingly fine discretizations of the problem.
One then uses a splitting iteration like Jacobi or symmetric successive overrelaxation (SSOR) to smooth different frequencies of the error $e = x^\ast - x$ on different grids, where $x^\ast$ solves $Ax^\ast = b$ and $x$ is an approximation.
Due to the decreasing computational complexity the problem can usually be solved explicitly on the coarsest grid, e.g. by a factorization of the system matrix.
Between the grids, interpolation and restriction operations are used in order to transport the information back and forth.
This works for hierarchical grids, i.e. that all nodes of one grid are also contained in the next finer grid.
We achieve this by successively halving the grid spacing from $2h$ to $h$.
A typical choice for the restriction is the transposed interpolation given by $I_{h}^{2h} = \left(I_{2h}^{h}\right){'}$.

The outstanding feature of the multigrid method is that under certain circumstances it is possible to solve sparse, linear systems in optimal $\mathcal{O}(m)$ complexity, where $m$ is the number of discretization points on the finest grid.
Algorithm \ref{Algorithm:MG} shows the basic structure of one V-cycle that serves as a preconditioner in the CG method.
Here $g \in \lbrace 1, \dots, G \rbrace$ denotes the grid levels from coarse $g=1$ to fine $g=G$, which coincides with the original problem.
The main ingredients of Algorithm \ref{Algorithm:MG} are:
\begin{itemize}
	\item the system matrices $A_g$,
	\item interpolation matrices $I_{g-1}^{g}$ and restriction matrices $I_{g}^{g-1}$,
	\item smoothing iteration method,
	\item the coarse grid solver $A_1^{-1}$.
\end{itemize}
Note that the smoother (line 6 and 11 in Algorithm \ref{Algorithm:MG}) does not necessarily have to be convergent on its own.
Typical choices are Jacobi or SSOR iterations, for which the convergence depends on the spectral radius of the iteration matrix.
Convergence is not necessary for the smoothing property.
It is thus recommendable only to apply a small number of pre/post-smoothing steps $\nu_1$ and $\nu_2$, such that a diverging smoother does not affect the overall convergence.
This is further addressed in the discussion of Section \ref{sec:numerics}.
\begin{algorithm}[htb]
	\caption{Multigrid v-cycle}
	\label{Algorithm:MG}
	\SetKwProg{Fn}{}{}{end}
	\Fn{\emph{v\_cycle}($A_g,b,g,x$)}{
		\If{$g=1$}{
			$x \gets A_g^{-1}b$ 
		}
		\Else{
			$x \gets \text{smooth}(x,b,\nu_1)$ \\
			$r \gets b - Ax$ \\
			$r \gets I_{g}^{g-1}r$ \\
			$e \gets \text{v\_cycle}(A_{g-1},r,g-1,0)$ \\
			$x \gets x + I_{g-1}^{g}e$ \\
			$x \gets \text{smooth}(x,b,\nu_2)$
		}
	}
\end{algorithm}

It might be tempting to apply an algebraic (AMG) instead of a geometric multigrid preconditioner.
The attractivity stems from the fact that AMG typically works as a black-box solver on a given matrix without relying a geometric description of grid levels and transfer operators.
Yet, common AMG implementations explicitly require access to the matrix, which is prohibitively memory consuming for the tensor-product smoothing splines.

In the following we concentrate on a memory efficient realization of the matrix-vector product in line 7 of Algorithm \ref{Algorithm:MG}, the interpolation/restriction in lines 8 and 10 and the smoothing operation in lines 6 and 11.

\subsection{Memory Efficient Matrix Operations}


Matrix-free methods for the solution of the large-scale linear system (\ref{Equation:SmoothingSplineLinearSystem}) require the efficient computation of matrix-vector products of the occurring matrices.
For our particular application, we exploit their special structure, namely Kronecker and Khatri-Rao product structure.


\subsubsection{Kronecker Matrices}


For given matrices $A_p \in \R^{m_p \times n_p}$ we consider the Kronecker matrix
\begin{align*}
A \DEF \bigotimes\limits_{p=1}^{P} A_p \in \R^{m \times n}, \ m \DEF \prod_{p=1}^{P} m_p , \ n \DEF \prod_{p=1}^{P} n_p
\end{align*}
and aim for a matrix-free computation of matrix-vector products with $A$ by only accessing the Kronecker factors $A_1,\ldots,A_P$.
In the case of $m_p=n_p$ for all $p=1,\ldots,P$ an implementation is provided by \cite{Benoit2001} and we extend their idea to arbitrary factors by Algorithm \ref{Algorithm:MVPKronecker} to form a matrix-vector product with a Kronecker matrix only depending on the Kronecker factors.
\begin{algorithm}[htb]
	\DontPrintSemicolon
	\caption{Matrix-Vector Product with Kronecker Matrix}
	\label{Algorithm:MVPKronecker}
	\KwIn{$A_1,\ldots,A_P, x$} 
	\KwOut{$v_0 = (A_1 \otimes \ldots \otimes A_P)x$}
	$v_P \gets x$ \\
	\For{$p=P,\ldots,1$}{ 
		\For{ $s=1,\ldots,l_p$}{
			$\bar{v}_{p,s} \gets v_p[(s-1)n_pr_p+1 : s n_pr_p]$ \\
			$w_{p,s} \gets 0$ \\
			\For{$t=1,\ldots,r_p$}{
				$z_{p,s,t} \gets \bar{v}_{p,s}[t,t+r_p,\ldots,t+(n_p-1)r_p]$ \\
				$\bar{z}_{p,s,t} \gets A_p z_{p,s,t}$ \\
				$w_{p,s}[t,t+r_p, \ldots, t+(m_p-1)r_p] \gets \bar{z}_{p,s,t}$
			}
			$v_{p-1}[(s-1)m_pr_p+1 : sm_pr_p ] \gets w_{p,s}$
		}
	}
\end{algorithm}


\subsubsection{Khatri-Rao Matrices}


For given matrices $A_p \in \R^{m_p \times n}$ we consider the Khatri-Rao matrix
\begin{align*}
A \DEF \bigodot\limits_{p=1}^{P} A_p \in \R^{m \times n}, \ m \DEF \prod_{p=1}^{P} m_p  
\end{align*}
and aim for a matrix-free computation of matrix-vector products with $A$ and $A{'}$ by only accessing the Khatri-Rao factors $A_1,\ldots,A_P$.
By definition of the Khatri-Rao product it holds
\begin{gather*}
A = \bigodot\limits_{p=1}^{P} A_p  = \begin{bmatrix} \bigotimes\limits_{p=1}^{P} A_p[\cdot,1] , \ldots , \bigotimes\limits_{p=1}^{P} A_p[\cdot,n] \end{bmatrix}
\end{gather*}
such that
\begin{align*}
Ax = \sum\limits_{i=1}^{n}x[i] v_i , \ v_i \DEF \bigotimes\limits_{p=1}^{P} A_{p}[\cdot,i] \in \R^{m} 
\end{align*}
for all $x \in \R^n$.
This yields Algorithm \ref{Algorithm:MVPKhatriRao} to form a memory-efficient matrix-vector product with a Khatri-Rao matrix only requiring the Khatri-Rao factors.
\begin{algorithm}[htb]
	\caption{Matrix-Vector Product with Khatri-Rao Matrix}
	\label{Algorithm:MVPKhatriRao}
	\KwIn{$A_1,\ldots,A_P,x$}
	\KwOut{$w = (A_1 \odot \ldots \odot A_P)x$}
	$w \gets 0$ \\
	\For{$i=1,\ldots,n$}{
		$v \gets A_1[\cdot,i] \otimes \ldots \otimes A_P[\cdot,i]$ \\
		$w \gets w+x[i]v$
	}
\end{algorithm}

Similarly, it holds
\begin{align*}
A{'}y = \begin{pmatrix} v_1{'} y  \\ \vdots  \\ v_n{'} y \end{pmatrix}
\end{align*}
for all $y \in \R^m$.
We can thus formulate Algorithm \ref{Algorithm:MVPKhatriRaoTransposed} to form a matrix-vector product with a transposed Khatri-Rao matrix by only accessing its Khatri-Rao factors.
\begin{algorithm}[htb]
	\caption{Matrix-Vector Product with Transposed Khatri-Rao Matrix}
	\label{Algorithm:MVPKhatriRaoTransposed}
	\KwIn{$A_1,\ldots,A_P,y$}
	\KwOut{$w \DEF (A_1 \odot \ldots \odot A_P)^{'}y$}
	$w \gets 0$ \\
	\For{$i=1,\ldots,n$}{
		$v \gets A_1[\cdot,i] \otimes \ldots \otimes A_P[\cdot,i]$ \\
		$w[i] \gets v{'}y$ 
	}
\end{algorithm}

For the implementation of the Jacobi smoother in the multigrid algorithm we need access to the diagonal of the matrix $AA^{'} \in \R^{m \times m}$ without computing the entire matrix but only accessing the Khatri-Rao factors $A_1,\ldots,A_P$.
For $j=1,\ldots,m$ the $j$-th diagonal element of $AA{'}$ is given by $e_j{'} AA{'} e_j = \Vert A{'} e_j \Vert_2^2$, where $e_j$ denotes the $j$-th unit vector, and it holds
\begin{align*}
\Vert A{'} e_j \Vert_2^2 = \left\Vert \begin{pmatrix} v_1{'} e_j  \\ \vdots  \\ v_n{'} e_j \end{pmatrix} \right\Vert_2^2 = \left\Vert \begin{pmatrix} v_1[j]  \\ \vdots  \\ v_n[j] \end{pmatrix} \right\Vert_2^2 = \sum\limits_{i=1}^{n} v_i[j]^2 .
\end{align*}
This yields Algorithm \ref{Algorithm:DiagonalKahtriRao} to extract the diagonal of $AA{'}$ by only accessing its factors.
\begin{algorithm}[htb]
	\caption{Diagonal of Product of Khatri-Rao Matrix and its Transposed}
	\label{Algorithm:DiagonalKahtriRao}
	\KwIn{$A_1,\ldots,A_P$}
	\KwOut{$d \DEF diag(AA{'})$, where $A \DEF A_1 \odot \ldots \odot A_P$}
	$d \gets 0$ \\
	\For{$i=1,\ldots,n$}{
		$v \gets A_1[\cdot,i] \otimes \ldots \otimes A_P[\cdot,i]$ \\
		\For{$j=1,\ldots,m$}{
			$d[j] \gets d[j] + v[j]^2$
		}
	}
\end{algorithm}


\subsection{MGCG Algorithm for Spline Smoothing}


We are now prepared to apply the multigrid V-cycle of Algorithm \ref{Algorithm:MG} as a preconditioner in the CG method to solve the linear system (\ref{Equation:SmoothingSplineLinearSystem}).

\subsubsection{Hierarchy}

For the multigrid method we require hierarchical grids denoted by $g=1,\ldots,G$, where $g=1$ is the coarsest and $g=G$ the finest.
Since we do not assume initial knowledge on the data we propose to base the underlying spline space on $m_p = 2^g-1$ equidistant knots in each space dimension $p=1,\ldots,P$.
For $g=1,\ldots,G$ we then define the coefficient matrix
\begin{align*}
	A_{g} \DEF \Phi{'}_g \Phi_g + \lambda \Lambda_g.
\end{align*}
With
\begin{align*}
	K_g \DEF dim(\mathcal{S}_q(\K_{g})) = \prod\limits_{p=1}^{P} (2^g + q_p) = (2^g + 3)^P
\end{align*}
we denote the dimension of the spline space $\mathcal{S}_q(\K_{g})$ and hence the size of the coefficient matrix $A_g \in \R^{K_g \times K_g}$ on grid level $g$.
Thus, we obtain a hierarchy of linear systems for which we can use the subdivision properties given in Section \ref{sec:regression}.

\subsubsection{Smoothing Iteration}

To apply the Jacobi method we additionally require the diagonal of the coefficient matrix, which is a vector of length $K_g$.
In principle, this is not prohibitively memory consuming yet, since the coefficient matrix is not explicitly accessible, we cannot simply extract its diagonal.
However, Algorithm \ref{Algorithm:DiagonalKahtriRao} allows the memory efficient computation of the diagonal of $\Phi{'}_g \Phi_g$, $g=1,\ldots,G$, and the diagonal of each matrix $\Psi_{r,g}$ is directly given by the diagonal property of Kronecker matrices.
In contrast to the Jacobi smoother the SSOR method additionally requires explicit access to all elements in one of the triangular parts of the coefficient matrices.
It would in principal be possible to compute the desired elements entry-wise within each iteration.
Yet, this would be computationally very expensive, since these elements have to be computed repeatedly in each iteration and for each grid level.
Therefore, we utilize the Jacobi method \ref{Algorithm:MemoryEfficientJacobi} as smoothing iteration.
\begin{algorithm}[htb]
	\DontPrintSemicolon
	\caption{Memory Efficient Jacobi Iteration for Spline Smoothing}
	\label{Algorithm:MemoryEfficientJacobi}
	\SetKwProg{Fn}{}{}{end}
	\Fn{JAC($\alpha,b,g,\nu$)}{
		$D_{\text{inv}} \gets 1 / diag( \Phi_g{'}\Phi_g + \lambda \Lambda_g)$ \tcp*[f]{Algorithm \ref{Algorithm:DiagonalKahtriRao}} \\
		\For{$j=1,\ldots,\nu$}{
			$r \gets b - \left( \Phi_g{'}\Phi_g + \lambda \Lambda_g \right)\alpha$ \tcp*[f]{Algorithm \ref{Algorithm:MVPKhatriRaoTransposed}, 	\ref{Algorithm:MVPKhatriRao}, and \ref{Algorithm:MVPKronecker}} \\
			$\alpha \gets \alpha + \omega D_{\text{inv}} r $
		}
	}
\end{algorithm}
We also test a SSOR smoother in the multigrid algorithm in Section \ref{sec:numerics}, yet we apply this only to low dimensional tests, where we can explicitly assemble and keep the coefficient matrix in memory.

\subsubsection{Grid Transfer}

For the transfer of a data vector from grid $g$ to grid $g+1$ we require a prolongation matrix $I_{g}^{g+1} \in \R^{K_{g+1} \times K_{g}}$.
Due to the subdivision formula and the tensor product nature of the splines we can use (\ref{Equation:SubdivisionMatrixMulti}) for this purpose.
For the restriction matrix $I_{g+1}^{g} \in \R^{K_{g} \times K_{g+1}}$ we choose $I_{g+1}^{g} \DEF \left( I_{g}^{g+1} \right){'}$, which yields the Garlerkin property
\begin{align*}
	A_{g} = I_{g+1}^{g} A_{g+1} I_{g}^{g+1} .
\end{align*}
The restriction and prolongation matrices do not fit into memory as well, but to apply the V-cycle Algorithm \ref{Algorithm:MG} only matrix-vector products with $I_{g}^{g+1}$ and $I_{g+1}^{g}$ are required.
Since
\begin{align*}
	I_{g}^{g+1} = \bigotimes_{p=1}^{P} I_{g,p}^{g+1,p} \ \text{ and } \ I_{g+1}^{g} \DEF \left( I_{g}^{g+1} \right){'} = \left( \bigotimes_{p=1}^{P} I_{g,p}^{g+1,p} \right){'} =  \bigotimes_{p=1}^{P} \left(I_{g,p}^{g+1,p} \right){'} ,
\end{align*}
these product can be memory efficiently computed by Algorithm \ref{Algorithm:MVPKronecker}.

\subsubsection{Coarse Grid Solver}

On the coarsest grid $g=1$ the V-cycle algorithm requires the exact solution of a linear system with coefficient matrix $A_{1} \in \R^{K_1 \times K_1}$.
Since $K_1 \ll K_G$ we assume an explicitly assembled coarse grid coefficient matrix in a sparse matrix format.
This allows a factorization of $A_{1}$, that can precomputed and stored in memory.
Keeping the factorization in memory is important, since each call of the multigrid preconditioner requires a solution of a linear system given by $A_{1}$.
Due to the symmetry of the matrix we apply a sparse Cholesky factorization.
Yet, in higher dimensional spaces even the coarse grid operators might be prohibitively memory consuming.
Recall that the number of variables $K_1$ grows exponentionally with the space dimension $P$.
We then apply the matrix-free CG algorithm also as a coarse grid solver.
Note that the overall scalability is not affected since the computational costs of the coarse grid solver are fixed, even when the fine grid $G$ is further subdivided as $G \gets G+1$.

\subsubsection{V-Cycle}

Putting everything together we obtain Algorithm \ref{Algorithm:MemoryEfficientMG}, which performs one V-cycle of the multigrid method for the large-scale linear system (\ref{Equation:SmoothingSplineLinearSystem}) with negligible memory requirement.
\begin{algorithm}[htb]
	\DontPrintSemicolon
	\caption{Memory Efficient V-Cycle for Spline Smoothing}
	\label{Algorithm:MemoryEfficientMG}
	\SetKwProg{Fn}{}{}{end}
	\Fn{v\_cycle($\alpha,b,g,\nu$)}{
		\If{$g=1$}{
			$\alpha \gets \left( \Phi_1{'}\Phi_1 + \lambda \Lambda_1 \right)^{-1}b$  
		}
		\Else{
			$\alpha \gets \text{JAC}(\alpha,b,g,\nu_1)$  \tcp*{Algorithm \ref{Algorithm:MemoryEfficientJacobi}}
			$r \gets \left( \Phi_g{'}\Phi_g + \lambda \Lambda_g \right)\alpha - b$  \tcp*{Algorithm \ref{Algorithm:MVPKhatriRaoTransposed}, 	\ref{Algorithm:MVPKhatriRao}, and \ref{Algorithm:MVPKronecker}}
			$r \gets I_{g}^{g-1}r$  \tcp*{Algorithm \ref{Algorithm:MVPKronecker}}
			$e \gets \text{v\_cycle}(0,r,g-1,\nu)$ \tcp*{Algorithm \ref{Algorithm:MemoryEfficientMG}}
			$\alpha \gets \alpha - I_{g-1}^{g}e$  \tcp*{Algorithm \ref{Algorithm:MVPKronecker}}
			$\alpha \gets \text{JAC}(\alpha,b,g,\nu_2)$  \tcp*{Algorithm \ref{Algorithm:MemoryEfficientJacobi}}
		}
	}
\end{algorithm}
The V-cycle can be interpreted as linear iteration with iteration matrix
\begin{align} \label{MGIterationMatrix}
	\begin{split}
		C_{\text{MG},G} &= C_{\text{smooth}}^{\nu_2} \left( I_{n_G} - I_{G-1}^{G} \left( I_{n_{G-1}} - C_{\text{MG},G-1} \right) A_{G-1}^{-1} I_{G}^{G-1} A_G \right) C_{\text{smooth}}^{\nu_1} , \\
		C_{\text{MG},1} &= 0 ,
	\end{split}
\end{align}
where $C_{\text{smooth}}$ denotes the iteration matrix of the utilized smoothing iteration \cite[cf.][p. 446]{Saad2003}.
Note that the iteration matrix is only for theoretical investigations and is never assembled in practical applications.

\subsubsection{MGCG Method}

Finally, applying the multigrid V-cycle as preconditioner for the CG method yields Algorithm \ref{Algorithm:MemoryEfficientMGCG} as memory efficient multigrid preconditioned conjugated gradient (MGCG) method to solve the large-scale linear system (\ref{Equation:SmoothingSplineLinearSystem}).
\begin{algorithm}[htb]
	\DontPrintSemicolon
	\caption{Memory Efficient MGCG Method for Spline Smoothing}
	\label{Algorithm:MemoryEfficientMGCG}
	$r \gets \Phi_G{'}y$ \tcp*{Algorithm \ref{Algorithm:MVPKhatriRao}}
	$p \gets z \gets \text{ v\_cycle}(0,r,G,\nu)$ \tcp*{Algorithm \ref{Algorithm:MemoryEfficientMG}}
	\While{stopping criterion not reached}{
		$v \gets  \left( \Phi_G{'}\Phi_G + \lambda \Lambda_G \right) p$ \tcp*{Algorithm \ref{Algorithm:MVPKhatriRaoTransposed}, \ref{Algorithm:MVPKhatriRao}, and \ref{Algorithm:MVPKronecker}}
		$w \gets \Vert r \Vert_2^2 / p{'}v$ \\
		$\alpha \gets \alpha + w p$ \\
		$\tilde{r} \gets r$ \\
		$r \gets r - w v$ \\
		$\tilde{z} \gets z$ \\
		$z \gets \text{v\_cycle}(0,r,G,\nu)$ \tcp*{Algorithm \ref{Algorithm:MemoryEfficientMG}}
		$p \gets z +  (r{'}z / \tilde{r}{'}\tilde{z}) p$
	}
\end{algorithm}


\section{Numerical Results}
\label{sec:numerics}

\begin{figure}[htb]
	\centering
	\begin{minipage}{0.45\textwidth}
		\includegraphics[width=\textwidth]{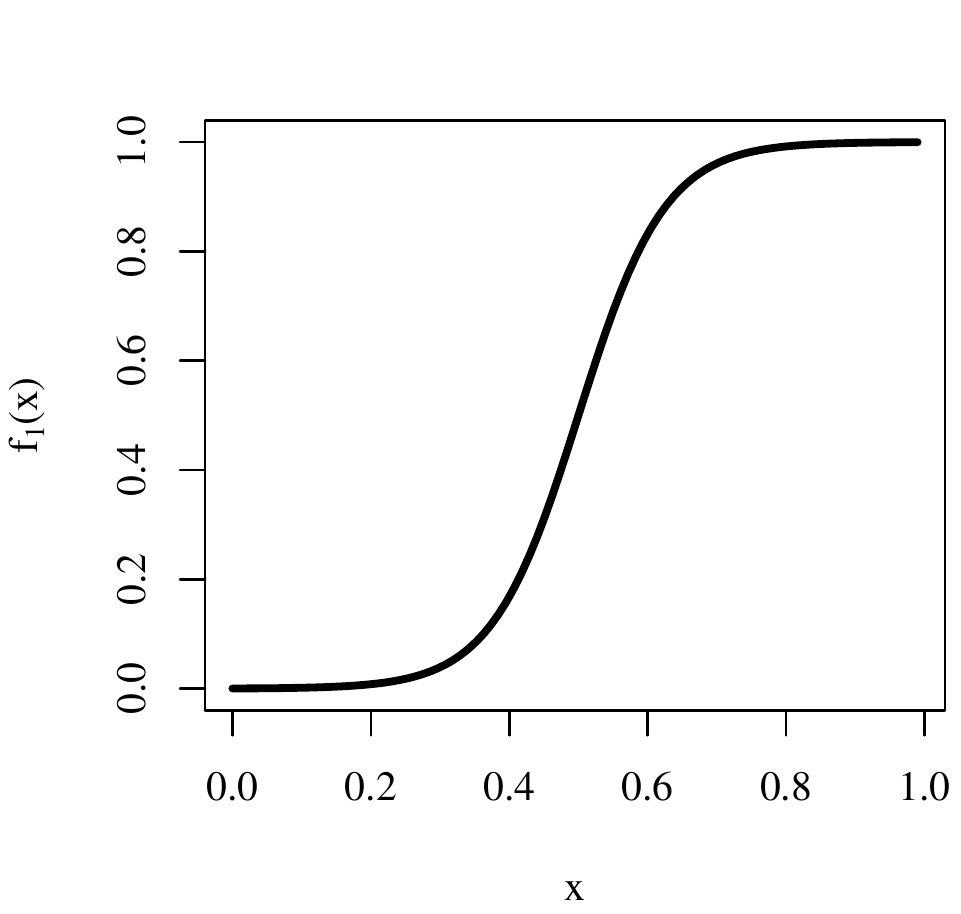}
	\end{minipage}
	\hspace{0.5cm}
	\begin{minipage}{0.45\textwidth}
		\includegraphics[width=\textwidth]{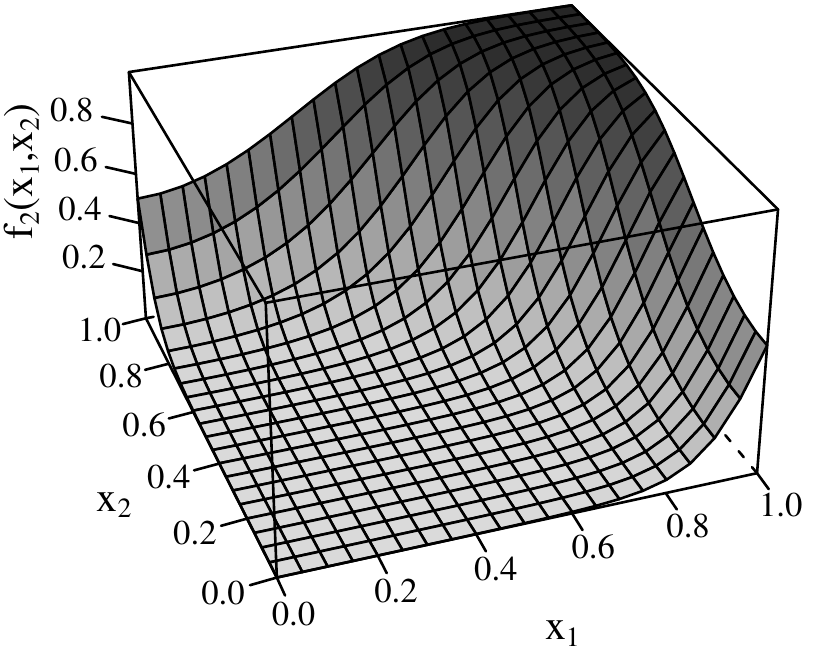}
	\end{minipage}
	\caption{Plot of the (undisturbed) sigmoid test function for $P=1$ and $P=2$.}
	\label{Figure:test_case}
\end{figure}

In this section we demonstrate the performance of Algorithm \ref{Algorithm:MemoryEfficientMGCG} in different numerical test cases.
On the one hand, we consider various spatial dimensions $P=1,\dots, 4$ in order to investigate the computational complexity of the algorithm.
On the other hand, we successively apply uniform grid refinements for a maximum grid level of $G=4, \dots, 7$ in each dimension in order to inspect the scalability of the multigrid preconditioner.
The underlying codes are programmed within R (version 3.4.4).
The algorithmic building blocks, which are critical to computational performance, are accelerated by using the RCPP extension library \cite[cf.][]{Eddelbuettel2011,Eddelbuettel2013} and programmed in {C++}.
These are in particular the matrix-vector products with the coefficient matrix (Algorithm \ref{Algorithm:MVPKronecker}, \ref{Algorithm:MVPKhatriRao}, \ref{Algorithm:MVPKhatriRaoTransposed}), interpolation and restriction (Algorithm \ref{Algorithm:MVPKronecker}), and the Jacobi smoother (Algorithm \ref{Algorithm:MemoryEfficientJacobi}).
Throughout this section, we consider the following test data set
\begin{align*}
	\lbrace (x_i,y_i) \in \R^P \times \R : i=1,\ldots,100{.}000 \rbrace , \ P=1,\ldots,4 , 
\end{align*}
obtained by uniformly random sampling of the normalized, multivariable sigmoid function
\begin{align*}
	f_P \colon [0,1]^P \rightarrow [0,1] , \ x \mapsto \frac{1}{1-\exp\left( -16 \left( \Vert x \Vert_2^2 P^{-1} - 0.5 \right) \right) }
\end{align*}
enriched by normally distributed noise $\varepsilon \sim \mathcal{N}(0;0.1^2)$, that is
\begin{align*}
	y_i \DEF f_P(x_i) + \varepsilon_i .
\end{align*}
For $P=1$ and $P=2$ the related sigmoid functions are graphed in Figure \ref{Figure:test_case}.
We consider these functions since they are irrational functions that show a similar behavior for varying covariate dimensions $P$.
However, since the performance of the solution algorithms is in focus, the exact form of the generating function is of minor importance.

\begin{figure}[htb]
	\centering
	\includegraphics[width=0.8\textwidth]{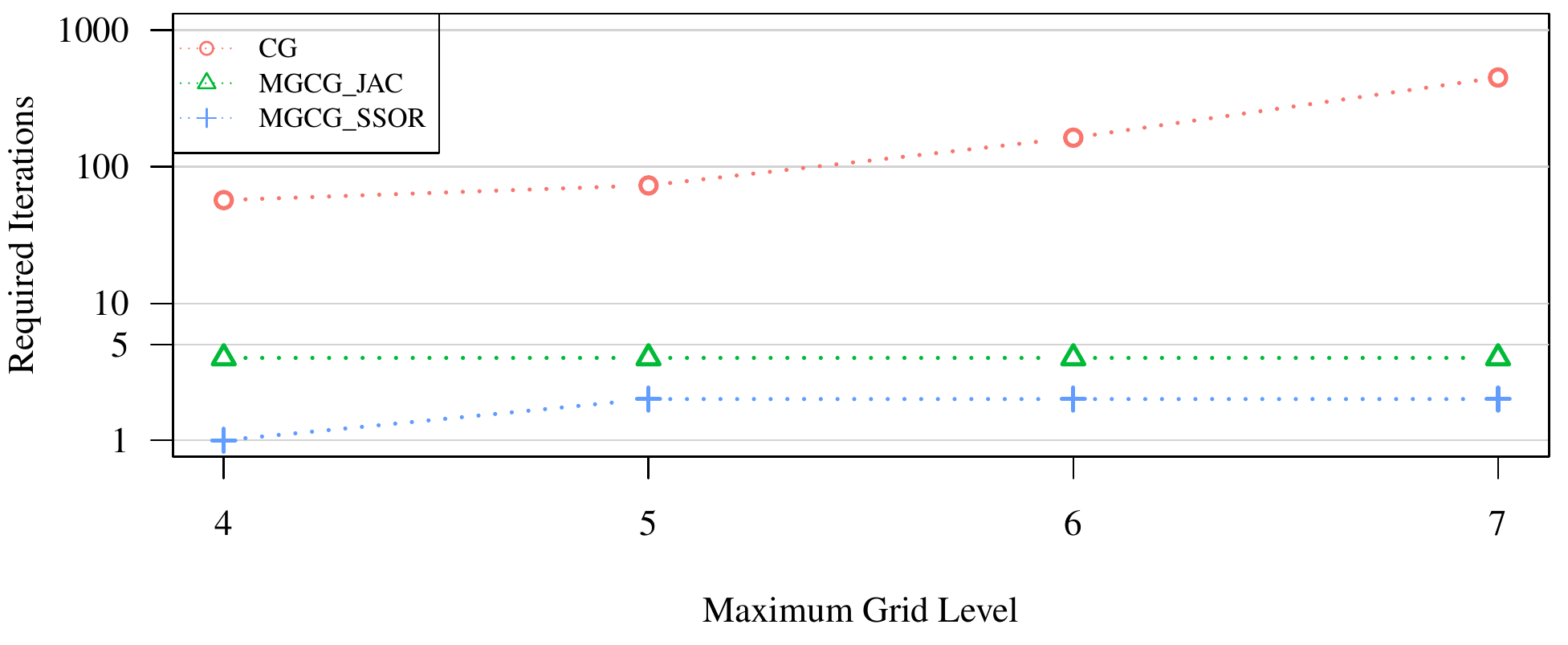}
	\caption{Number of preconditioned and unpreconditioned conjugate gradient iterations under uniform grid refinements in $P=2$ dimensions.}
	\label{Figure:grid_refinements}
\end{figure}

In a first test case, we fix the spatial dimension to $P=2$ and successively refine the B-spline basis on $G=1,\dots,7$ grids.
We then consider four different tests.
In each of them, $G \in \lbrace 4,\dots,7\rbrace$ is chosen to be the maximum grid level for the V-cycle of the multigrid preconditioner.
The coarsest grid $g=1$ is used for the coarse grid solver in each case.
This leads to problem dimensions of $K_1 = (2^1+3)^2=25$ on the coarse grid and $K_G = (2^G+3)^2$ on the finest grid.
Note that the unpreconditioned CG method simply uses the discretization matrix $A_G$ on the $G$-th grid.
The results are shown in Figure \ref{Figure:grid_refinements}, where a logarithmic scale is used.
We observe that the number of unpreconditioned conjugate gradient iterations significantly increases under grid refinements, whereas for both, the Jacobi and SSOR preconditioned MGCG algorithms, the number of iterations is almost constant (i.e. 1-2 for SSOR and 4 for Jacobi).
This result illustrates that the multigrid preconditioner enables a scalable solver for the regularized least squares problem (\ref{Equation:SmoothingSplineLinearSystem}) determining the smoothing spline.
Clearly, the computational times are increasing since we run only a single core code.
A standard approach here is to distribute the matrix between an increasing number of processors in a cluster computer and one obtains almost constant running times in the sense of weak scalability.
Note that this is not achievable by the plain CG solver due to the increasing number of iterations.

The next test under consideration is the comparison of computational complexity with respect to the spatial dimension.
Table \ref{Table:computational_time} shows the number of required iterations of the utilized methods in the dimensions $P=1,\dots,4$ together with the CPU time in seconds.
The underlying hardware is one core of a Xeon W-2155 processor.
One can observe that each of the preconditioned iterations is more expensive with respect to computational time.
Yet, the MGCG is significantly faster due to the small number of iterations. 
Note that the computational time of the SSOR multigrid preconditioned iteration is not comparable since it does not rely on the matrix-free approach.
This is only for the comparison of MGCG iterations.
Since the SSOR smoother needs access to all matrix entries, we explicitly assemble the coefficient matrix $\left( \Phi{'} \Phi + \lambda \Lambda \right)$ in R's sparse format and use it for the algorithm.
\begin{table}[htb]
	\centering
	\begin{tabular}{r | r r | r r | r r}
		& \multicolumn{2}{c|}{CG} & \multicolumn{2}{c|}{MGCG\_JAC} & \multicolumn{2}{c}{MGCG\_SSOR} \\ \hline \raisebox{2.5ex}{}
		P=1 & 22 & (1.05) & 2 & (0.91) & 1 & (<0.01) \\
		P=2 & 73 & (4.41) & 4 & (4.59) & 2 & (0.37) \\
		P=3 & 367& (103.31) & 14& (79.9) & \multicolumn{2}{c}{-} \\
		P=4 & 747& (2745.98) & 19 & (1927.60) & \multicolumn{2}{c}{-}	
	\end{tabular}
	\caption{Required number of iterations of the considered methods for various dimensions with $G=5$ grids and computational time in seconds.}
	\label{Table:computational_time}
\end{table}

\begin{figure}[htb]
	\centering
	\includegraphics[width=1.0\textwidth]{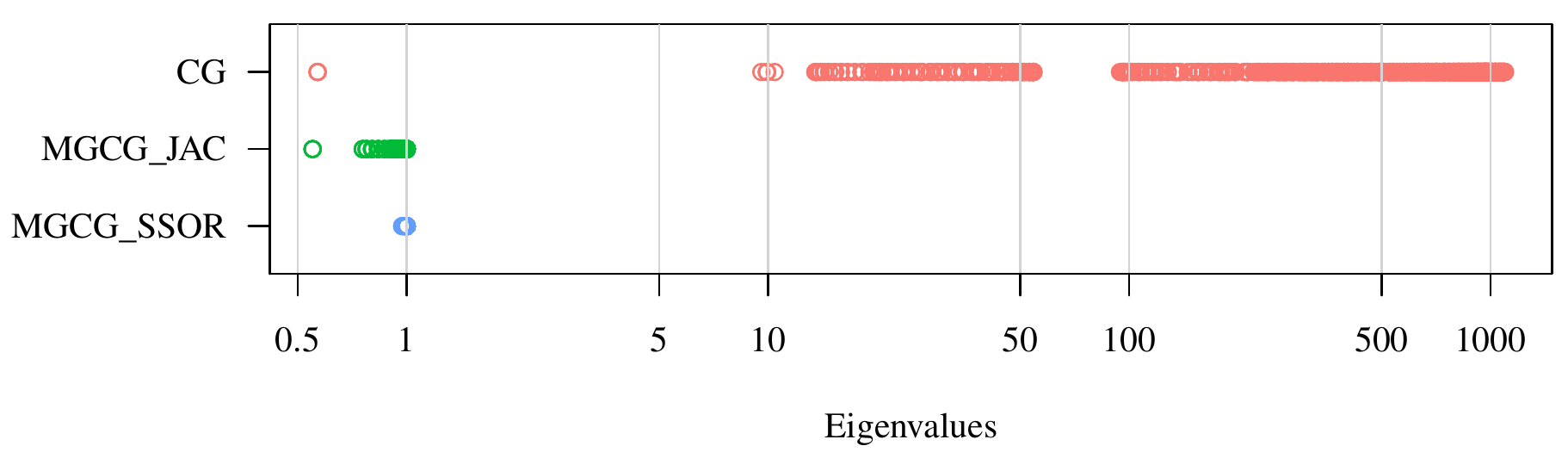}
	\caption{Eigenvalues of the (un)preconditioned coefficient matrices for $G=5$ grids and $P=2$ dimensions.}
	\label{Figure:eigenvalues}
\end{figure}

Crucial for the convergence speed of the CG method in the unpreconditioned as well as in the preconditioned case is the condition of the respective coefficient matrix, that is
\begin{align*}
\Phi_G{'}\Phi_G + \lambda \Lambda_G
\end{align*}
for the pure CG method and 
\begin{align*}
C_{\text{MG},G} \cdot \left( \Phi_G{'}\Phi_G + \lambda \Lambda_G \right)
\end{align*}
for the MGCG method.
An explicit form of the iteration matrix of the multigrid method is given in (\ref{MGIterationMatrix}), i.e. one call of Algorithm \ref{Algorithm:MG} is equivalent to a matrix-vector product with $C_{\text{MG},G}$.
Note that for the the Jacobi smoother and for the SSOR smoother different preconditioning matrices $C_{\text{MG},G,\text{JAC}}$ and $C_{\text{MG},G,\text{SSOR}}$ arise.
Figure \ref{Figure:eigenvalues} visualizes the distribution of the eigenvalues of the corresponding coefficient matrices for the $P=2$ and $G=5$ case.
Additionally, Table \ref{Table:condition} shows the condition number of the (un)preconditioned iteration matrix for the CG algorithm.
Note that we only compute the eigenvalues in the $P=2$ test case due to the computational complexity, since we have to assemble a matrix representation of the preconditiond system matrix.
That is, we have to run the V-cycle on each of the $K$-dimensional unit vectors to obtain a matrix on which we then perform an eigenvalue decomposition.
Here it can be seen that the multigrid preconditioner pushes the eigenvalues towards 1, which explains the significant lower number of required MGCG iterations.
Although it seems that the multigrid preconditioner with SSOR smoother outperforms the Jacobi smoother, it is prohibitively expensive with respect to memory requirement.
For the SSOR iteration the entire upper triangular part of the system matrix is required, which is not efficiently accessible in our matrix-free approach.
In $P=3$ dimensions, the MGCG with SSOR smoother for the considered test problem requires approximately 30 GB of RAM, which is at the limit of the utilized computer system.
In contrast, the matrix-free methods require approximately 16 MB (CG) and 78 MB (MGCG\_JAC) of RAM.
Due to the fact, that in the matrix free approach we have computational redundancy, the computational times of the full approach are significantly better.
Yet, this is bounded to the low dimensional case because of the exponential growth of required memory.
\begin{table}[htb]
	\centering
	\begin{tabular}{c|c|c|c}
		& CG & MGCG\_JAC & MGCG\_SSOR \\ \hline
		condition & $1933.27$ & $1.82$ & $1.03$
	\end{tabular}
	\caption{Condition number of the coefficient matrices for $G=5$ grids and $P=2$ dimensions.}
	\label{Table:condition}
\end{table}

\begin{figure}[htb]
	\centering
	\begin{minipage}{0.49\textwidth}
		\includegraphics[width=\textwidth]{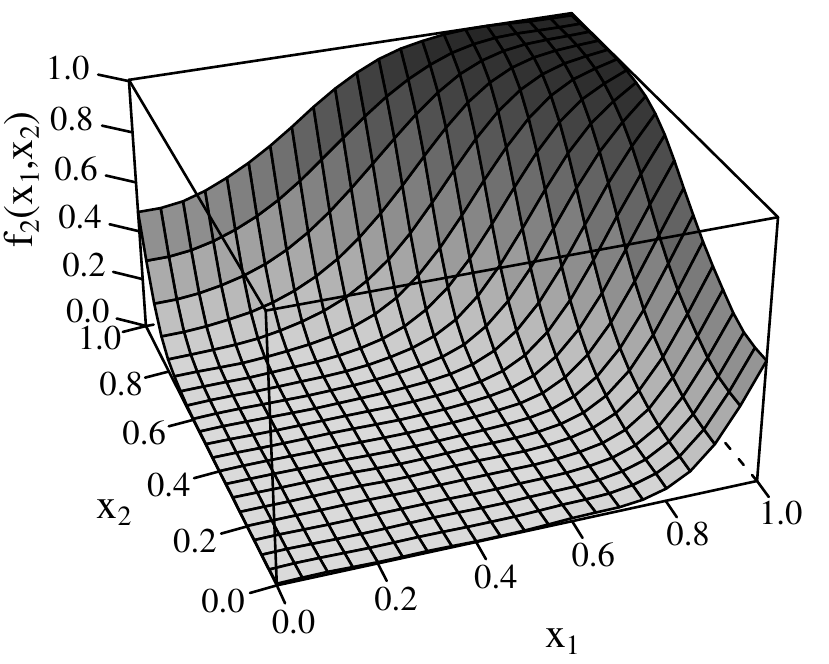}
	\end{minipage}
	\begin{minipage}{0.49\textwidth}
		\includegraphics[width=\textwidth]{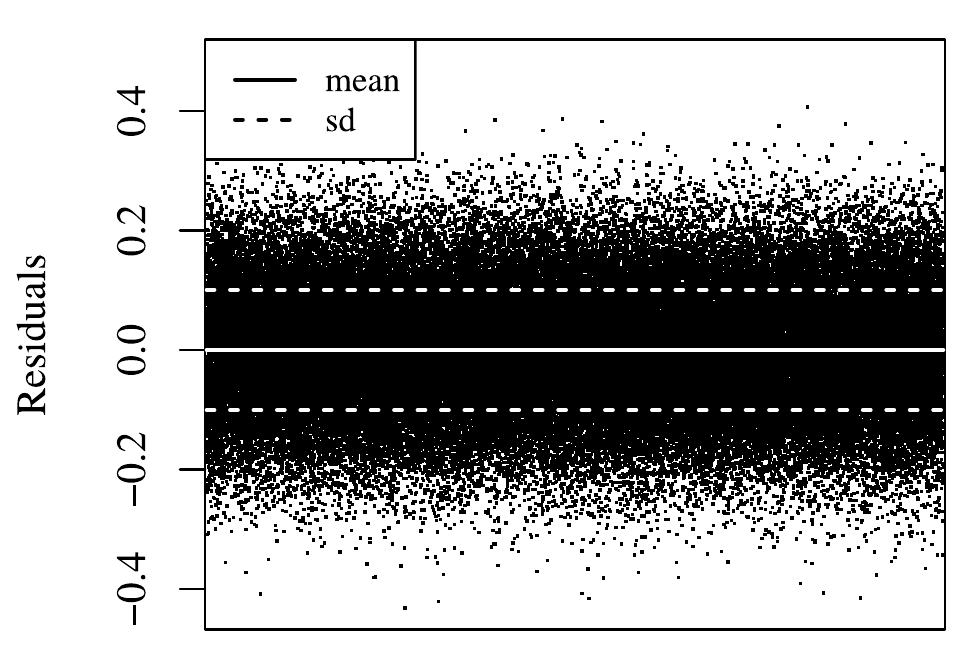}
	\end{minipage}
	\caption{Smoothing spline approximation (left) and related residuals of the approximation (right) for $G=5$ grids and $P=2$ dimensions.}
	\label{Figure:results}
\end{figure}

Finally, Figure \ref{Figure:results} shows the results of the smoothing spline approximation in two dimensions on the left-hand side and its corresponding residuals to the $100{.}000$ noisy data points on the right-hand side.
From this we can see, that the resulting smoothing spline recovers the underlying test function with adequate precision.

\section{Conclusions}
\label{sec:conclusion}

In this paper, we have presented a memory efficient algorithm in order to determine a smoothing spline in increasing spatial dimensions.
An important feature of our approach is the possibility to handle also scattered data in contrast to the already existing methods for gridded data.
The main challenge is to deal with memory and computational complexity originating from the large-scale linear system arising from spline smoothing with scattered data sets.
In order to overcome the issue of memory requirements, we have initially implemented a matrix-free conjugate gradient method, which comes along without assembling and storing the occurring matrices, but relies solely on matrix-vector products.
For this purpose, we especially have exploited the inherent tensor product structure of the multivariable spline functions.
Moreover, we address the issue of computational complexity by applying a geometric multigrid preconditioner to the CG algorithm.
This enables almost constant iteration numbers for fixed dimensions even under arbitrary grid refinements, which provides an important building block for algorithmic scalability.
In a representative numerical case study, we show the applicability and performance of the proposed method.


\section*{Acknowledgment}
This work has been partly supported by the German Research Foundation (DFG) within the research training group ALOP (GRK 2126).




\end{document}